\documentstyle{amsart}

\newtheorem{thm}{Theorem}[section]
\newtheorem{prop}[thm]{Proposition}
\newtheorem{lem}[thm]{Lemma}

\theoremstyle{definition}
\newtheorem{defn}[thm]{Definition}
\newtheorem{remark}[thm]{Remark}

\newtheorem{reduction}[thm]{Reduction}

\numberwithin{equation}{section}

\newcommand{\Spec}{\operatorname{Spec}}

\newcommand{\Hom}{\operatorname{Hom}}

\renewcommand{\Im}{\operatorname{Im}}

\newcommand{\bbA}{{\Bbb A}}

\newcommand{\bbC}{{\Bbb C}}
\newcommand{\bbZ}{{\Bbb Z}}

\newcommand{\bbQ}{{\Bbb Q}}

\newcommand{\bbP}{{\Bbb P}}

\newcommand{\cO}{{\cal O}}
\newcommand{\cM}{{\cal M}}

\newcommand{\cI}{{\cal I}}

\newcommand{\fm}{{\frak m}}

\newcommand{\kbar}{\overline k}

\newcommand{\isomo}{\overset{\sim}{=}}

\newcommand{\lf}{\mathopen}

\let\r=\mathclose

\newcommand{\Stab}{\operatorname{Stab}}

\newcommand{\lra}{\longrightarrow}

\let\to=\longrightarrow

\newcommand{\paru}{Parusi\'nski}
\newcommand{\projC}{\bbP^{n-1}_\bbC}
\newcommand{\projQ}{\bbP^{n-1}_\bbQ}

\newcommand{\nvstab}{\operatorname{NStab}}

\begin{document}
\title[Parusi\'nski's Lemma]{Parusi\'nski's ``Key Lemma" via algebraic geometry}
\author[Z. REICHSTEIN and B. YOUSSIN,  10-16-99]
{Z. Reichstein and B. Youssin} 
\address{Department of Mathematics, Oregon State University,
Corvallis, OR 97331} 
\thanks{Z. Reichstein was partially supported by NSF grant DMS-9801675 and
(during his stay at MSRI) by NSF grant DMS-9701755.}
\email{zinovy@@math.orst.edu}
\address{Department of Mathematics and Computer Science,
University of the Negev, Be'er Sheva', Israel\hfill\break
\hbox{{\rm\it\hskip\parindent Current mailing address\/}}: 
Hashofar 26/3, Ma'ale Adumim, Israel}
\email{youssin@@math.bgu.ac.il}
\subjclass{Primary 14E15, 14F10, 14L30. 
Secondary: 16S35, 32B10, 58A40}

\begin{abstract}
The following ``Key Lemma'' plays an important role in Parusinski's 
work on the existence of Lipschitz stratifications in the class of
semianalytic sets:
For any positive integer $n$, there is a finite set of homogeneous
symmetric polynomials $W_1, \dots ,W_N$ in $Z[x_1,...,x_n]$
and a constant $M >0$ such that
\[ |dx_i/x_i| \le M \max_{j = 1, \dots, N} |dW_j/W_j| \; , \]
as densely defined functions on the tangent bundle of $C^n$.
We give a new algebro-geometric proof of this result.
\end{abstract}

\maketitle

\section{Introduction}

\paru's fundamental work on the existence 
of Lipschitz stratifications in the class of semianalytic sets
relies on the following result.

\begin{thm}  (\paru~\cite[pp. 202--203]{P}) 
\label{thm2}
For any positive integer\/ $n$, there is a finite set of homogeneous
symmetric polynomials $W_1,\dots,W_N \in\bbZ[x_1,\dots,x_n]$
and a constant $M >0$ such that
\begin{equation}
\label{eqn100}
\Bigl|\frac{dx_i}{x_i}(p,v)\Bigr|\le
M \max_{j=1,\dots,N}\Bigl|\frac{dW_j}{W_j}(p,v)\Bigr|
\end{equation}
for all $p\in\bbC^n$ and $v\in T_p\bbC^n$ for which both sides are defined.
Here, for any $P \in \bbC[x_1, \ldots, x_n]$ we view the meromorphic
differential form $\frac{dP}{P}$ on $\bbC^n$ as a densely 
defined function on the total space of the tangent bundle $T\bbC^n$.
\end{thm}

\paru\ refers to Theorem~\ref{thm2} as the ``Key Lemma"; the proof
of this result in~\cite[Section~6]{P} is quite difficult, 
being apparently the hardest part of~\cite{P}. 
The purpose of this paper is to show that Theorem~\ref{thm2}, 
in spite of its analytic appearance, 
has a natural proof in the framework of algebraic geometry.
Our argument is an application of the results of~\cite{ry} about 
group actions on algebraic varieties; these results, in turn, rely 
on canonical resolution of singularities.

We remark that \paru\ proves 
the inequality~\eqref{eqn100} under the additional assumption that 
$dV(p,v)=0$ if $V(p)=0$ for every $V$ belonging to finite set $\cal{V}$
of polynomials. Since this additional requirement does not affect
a dense Zariski open subset of $T\bbC^m$ (given by $V(p) \neq 0$ for
every $V \in \cal{V}$), it can be dropped. 
We also note that the statement 
of the Key Lemma in~\cite{P} only asserts the existence of polynomials 
$W_1, \dots, W_N$ with real coefficients; however, the construction 
of $W_1, \dots, W_N$ given there,
produces polynomials over $\bbZ$.  
Thus, while Theorem~\ref{thm2} 
appears to be stronger than the ``Key Lemma'' 
in~\cite{P}, the two are, in fact, equivalent.

We are grateful to R. Guralnick and S. Montgomery for informative 
conversations on the subject of group algebras, and to
D. Kazhdan for helpful comments on 
an earlier draft of this paper.

\section{Preliminaries}

\subsection*{Notational conventions}
All algebraic varieties considered in this paper, are assumed to be
irreducible and definied over a field $k$ of characteristic $0$. 
The base field $k$ is not assumed to be algebraically closed; 
the two cases of interest to us are $k=\bbQ$ and $k=\bbC$.  
By a point of a variety 
we shall always understand a closed point.
Given an embedding $k\subset\bbC$ and
a rational function $f$ on $X$, we shall denote the corresponding 
rational function on $X_\bbC=X\otimes_k\bbC$ by $f$ as well.

%
Throughout this paper $G$ will be a finite group.
A $G$-variety $X$ is a variety with a regular action of $G$,
$G\times X\to X$, where $G\times X$ is understood as the disjoint
union of $|G|$ copies of $X$.
We will always assume that the $G$-action is faithful, 
i.e., every nonidentity element of $G$ acts nontrivially.
By a morphism (respectively, rational map, birational isomorphism)
of $G$-varieties we shall mean a $G$-equivariant morphism (respectively,
rational map, birational isomorphism).

\subsection*{Stabilizers}
For a point $x$ in a $G$-variety $X$, we define its ``naive''
stabilizer $\nvstab(x)$ as the set of all $g\in G$ which preserve $x$.
If $k$ is not algebraically closed, the residue field $k(x)$ may be a
nontrivial finite extension of $k$, and $G$ may act on it
nontrivially. We define the ``honest'' stabilizer $\Stab(x)$ as the
set of all $g\in\nvstab(x)$ that act on $k(x)$ trivially;
cf.~\cite[Definition~0.4]{git}.
The subgroups $\nvstab(x)$ and $\Stab(x)$ of $G$ are
sometimes called the {\em decomposition group\/} 
and the {\em inertia group\/} respectively.

If $\kbar$ is the algebraic closure of $k$, then $x$ is represented by
a set of ``conjugate" points of the variety $X_{\kbar}=X\otimes_k\kbar$ 
(``geometric points'' of $X$), one for each embedding
$k(x)\hookrightarrow\kbar$;
$\Stab(x)$ fixes each of these points
while $\nvstab(x)$ permutes them.
As an example, consider the action of $G = \bbZ/2\bbZ$ on
the affine line $\bbA^1_\bbQ=\Spec\bbQ[t]$: the nontrivial
element of this group acts by $t \mapsto -t$. Here $\Stab(x)=\{1\}$
for any $x\in\bbA^1_\bbQ-\{0\}$. On the other hand,
$\nvstab(x)=G$
iff $x$ corresponds to the ideal in $\bbQ[t]$ generated by an
irreducible polynomial of the form $q(t^2)$ (e.g., $t^2 + 1$).

This phenomenon is entirely arithmetic; we are concerned with it here
because the symmetric polynomials $W_1, \dots, W_N$ in Theorem~\ref{thm2}
are asserted to have integer coefficients.  A reader who is only interested
in the existence of such polynomials in $\bbC[x_1, \dots, x_n]$
may skip the rest of this section and assume 
that ``naive'' stabilizers always coincide with ``honest'' 
ones in the sequel. 

\subsection*{Semi-linear representations and skew group rings}

Let $X$ be a $G$-variety and let $x \in X$.
The ``naive stabilizer'' $\nvstab(x)$ acts upon $T_x(X)^*=\fm_x/\fm_x^2$.
However, if $\nvstab(x)$ is strictly larger than $\Stab(x)$ 
then this action is not linear over $k(x)$
but rather ``semi-linear" in the following sense.

\begin{defn} \label{def.skew-lin} 
Suppose a finite group $H$ acts by automorphisms on a field $K$.
A {\em semi-linear representation} of $H$ over $K$ 
is a $K$-vector space $V$ 
with a $K^H$-linear action of $H$ on $V$ having the
property $g(\lambda v)=g(\lambda)g(v)$ 
for any $g\in H$, $\lambda\in K$ and $v\in V$.
\end{defn}

For the rest of this section we shall assume that $K$ is a field,
$K^*$ is the multiplicative group of $K$,
$H$ is a finite group acting on $K$ by automorphisms, and $H'$
is the kernel of this action. In the subsequent applications we 
will take $K = k(x)$, $H = \nvstab(x)$ and $H' = \Stab(x)$.

Recall that the skew group algebra $K*H$ is defined as the set 
of formal sums $\sum_{h\in H}a_hh$ (where $a_h\in K$), 
with componentwise addition and with multiplication given, distributively,
by $(a_1h_1)(a_2h_2) = a_1h_1(a_2) h_1h_2$. 
A semi-linear representation of $H$ is the same thing 
as a $(K*H)$-module. (All modules in this paper are
understood to be left modules.)

\begin{remark} \label{rem:triv}
Note that $V=K$ has a natural structure of
a  $(K*H)$-module.  This module contains a vector
$1\in K$ which is fixed by $H$.
\end{remark}

Recall that by Wedderburn's Theorem 
every semisimple ring $R$ is a direct product 
of simple rings, called {\em the simple components} of $R$,
see, e.g., \cite[Theorem~VIII.5.1]{bourbaki}. 

\begin{lem} \label{lem.skew-gr1}
$K*H$ is a semisimple $K^H$-algebra with at most $|H'|$ simple
components.
(Here $H'$ is the kernel of the $H$-action on $K$, as above.)
\end{lem}

\begin{pf}  Semisimplicity of $K*H$ is proved by
the same averaging argument as the usual Maschke's theorem; for
details see, e.g.,~\cite[Theorem 0.1 and Corollary~0.2]{montgomery}.

Denote the simple components of $K*H$ by $S_1, \dots, S_m$. Then 
$Z(K*H)=Z(S_1) \times \dots \times Z(S_m)$,
where $Z(A)$ denotes the center of $A$.
It is easy to see directly that $\dim_{K^H}(K*H) \le |H'|$.
Hence, $m \leq |H'|$, as claimed.
%
\end{pf}

The following proposition describes the particular kind of skew
group rings we shall encounter in the sequel.

\begin{prop} \label{prop:one-dim}
Suppose that $H'$ is an abelian group of exponent $e$, $K$ contains 
a primitive $e$th root of unity
and $\chi_1,...,\chi_m$ is a set of generators for the dual
group $(H')^*=\Hom(H',K^*)$. Assume further that for each $i$ 
there is a one-dimensinal semi-linear representation
$V_i$ of $H$ over $K$ such that $h'(v) = \chi_i(h')v$ for every
$v \in V_i$ and every $h' \in H'$.  Then:
\begin{itemize}
\item[(a)] For every $\chi\in(H')^*$, there exists a unique
semi-linear representation $V_\chi$ of $H$ such that 
$\dim_K(V_{\chi}) = 1$ and $h'(v)=\chi(h')v$ 
for every $v \in V_\chi$ and every $h' \in H'$.
\item[(b)] Every simple 
$(K*H)$-module is isomorphic to $V_{\chi}$ for some $\chi \in(H')^*$.
\end{itemize}
\end{prop}

\begin{pf} 
Note that if $V_1$ and $V_2$ are semi-linear representations of $H$
over $K$ then so is $V_1 \otimes_K V_2$. Indeed,
\[ h(\lambda v_1)\otimes h(v_2)=h(\lambda)\cdot(h(v_1)\otimes h(v_2))=
h(v_1)\otimes h(\lambda v_2) \; . \]
To construct $V_\chi$,
write $\chi \in (H')^*$ as
$\chi=\chi_1^{l_1}\dots\chi_m^{l_m}$ 
for some nonnegative integers $l_1,\dots,l_m$, and set
$ V_\chi=
V_1^{\otimes l_1}\otimes_{K}\dots\otimes_{K}V_m^{\otimes l_m}$.
The subgroup $H'$ acts on $V_\chi$ by the character $\chi$, 
as desired.  As $\dim_K(V_\chi) = 1$,
$V_\chi$ is an simple  $(K*H)$-module. 
Note that the $(K*H)$-modules $V_{\chi}$ 
are pairwise nonisomorphic because $H'$ acts on them 
by different characters.

The isomorphism classes of simple $(K*H)$-modules are in 1---1
correspondence with the simple components of $K*H$; see
\cite[Proposition~VIII.5.11]{bourbaki}. Thus Lemma~\ref{lem.skew-gr1}
implies that $K*H$ has $\leq |H'|$ nonisomorphic simple modules.
On the other hand, we have constructed $|H'|$ nonisomorphic
simple modules $V_{\chi}$.  This proves (b) and 
the uniqueness of $V_\chi$ in (a).
\end{pf}

\begin{remark} \label{rem.H} 
One can show that, under the assumptions of 
Proposition~\ref{prop:one-dim},
$H$ is a semidirect product of $H'$ and $H/H'$, where the
the action of $H/H'$ on $H'$ is given by embedding $H'$ into
$(K^*)^m$ via $h' \longmapsto (\chi_1(h'), \dots, \chi_m(h'))$. 
\end{remark}

\section{Reduction to an algebro-geometric problem}
\label{sect2}

%
%
%

We begin by restating~\eqref{eqn100} as an inequality involving
densely defined functions on the tangent bundle of $\projC$ rather
than the tangent bundle of $\bbC^n$.
Since $\projC$ is compact in the metric topology,
this will allow us to pass from local to global estimates.

\begin{prop} \label{prop2a} 
Let $X$ be a projective $G$-variety over a field $k\subset\bbC$, and
let $f$ be
a rational function on $X$.  Then there exist $G$-invariant rational
functions $\beta_1, \dots, \beta_m$ on $X$ and a constant $K > 0$ such that 
\begin{equation} \label{eqn101} 
\Bigl| \frac{df}{f}(p, v) \Bigr|\le
K \max_{j=1,\dots,m}\Bigl|\frac{d\beta_j}{\beta_j}(p,v)\Bigr|
\end{equation}
for any $(p, v) \in T(X_\bbC)$ such that $p$ is a smooth point of
$X_\bbC$ and does not lie on
the divisors of $f, \beta_1, \dots, \beta_m$.
\end{prop}

\begin{reduction} \label{red2} {\em Proposition~\ref{prop2a} $\Longrightarrow$
Theorem~\ref{thm2}.}

\smallskip
Indeed, apply Proposition~\ref{prop2a} with $k = \bbQ$,
$X = \projQ$, $G=S_n$ and $f = x_1/s_1$; here $S_n$ acts on $\projQ$ by
permutations of the homogeneous coordinates $x_1, \dots, x_n$ and
$s_1 = x_1 + \dots + x_n$.
Write each $\beta_j$ as a quotient of two homogeneous polynomials (of 
the same degree, with integer coefficients) in $x_1, \dots, x_n$:
\[ \beta_1 = W_1/W_2\,,\, \dots\,,\, \beta_m = W_{2m-1}/W_{2m} \; . \]
We claim that the polynomials $W_1, \dots, W_{2m}, W_{2m+1} 
\stackrel{\text{def}}{=} s_1$ have the property asserted in Theorem~\ref{thm2}.
Indeed, since $df/f = dx_1/x_1 - ds_1/s_1$ and 
$d\beta_i/\beta_i = dW_{2i-1}/W_{2i-1} - dW_{2i}/W_{2i}$, inequality
\eqref{eqn101} translates into
\begin{multline*}
\Bigl| \left(\frac{dx_1}{x_1} - \frac{ds_1}{s_1}\right)(p, v)  
\Bigr| \le
K \max_{i=1,\dots,m}\Bigl| \left(\frac{dW_{2i-1}}{W_{2i-1}} - 
\frac{dW_{2i}}{W_{2i}}\right)(p, v) \Bigr| \\
 \le 2K\max_{j=1,\dots,2m}\Bigl| \frac{dW_j}{W_j}(p, v) \Bigr|\ .
\end{multline*}
Consequently, 
\begin{multline*} \Bigl| \frac{dx_1}{x_1} (p, v) \Bigr| \le 
2K\max_{j=1,\dots,2m}\Bigl| \frac{dW_j}{W_j}(p, v) \Bigr| + 
\Bigl| \frac{ds_1}{s_1}(p, v) \Bigr|
\\ \le (2K+1)\max_{j=1,\dots,2m+1}\Bigl| \frac{dW_j}{W_j}(p, v) \Bigr| \; . 
\end{multline*}
This means that \eqref{eqn100} holds for $i = 1$,
with $M = 2K+1$ and $N = 2m+1$.
By symmetry, \eqref{eqn100} holds for all $i$.
\qed
\end{reduction}

\begin{defn} \label{def.prop-star}
Let $X$ be a projective $G$-variety and $f$ be a rational function on $X$.
We shall say that the pair $(X, f)$ has property (*) if
there is a Zariski open covering $X=\bigcup_iU_i$ and, for each $i$,
rational functions $\beta_{i1}, \dots, \beta_{i,q_i} \in k(X)^G$ and
regular functions $\gamma_{i1},\dots,\gamma_{i,q_i}\in \cO_X(U_i)$
such that
\begin{equation} \label{e.star}
\frac{df}{f} = \gamma_{i1}\frac{d\beta_{i1}}{\beta_{i1}} + \dots
\gamma_{i,q_i}\frac{d\beta_{i,q_i}}{\beta_{i,q_i}}  \; .  
\end{equation}
In other words, the pair $(X, f)$ has property (*) if $df/f$ is a global 
section of the sheaf of differentials 
on $X$ generated over $\cO_X$
by $d\beta/\beta$, as $\beta$ ranges over some finite subset of $k(X)^G$
(or, equivalently, as $\beta$ ranges over all of $k(X)^G$).
\end{defn}

\begin{reduction} \label{red3}
{\em Proposition~\ref{prop2a} holds, assuming the pair
$(X,f)$ that appears there, has property (*).}

\smallskip
Indeed, the Zariski open cover 
$\bigcup_iU_i$ of $X$, as in Definition~\ref{def.prop-star},
gives rise to
a Zariski open cover $\bigcup_iU_{i,\bbC}$ of $X_{\bbC}$.
The functions $\gamma_{ij}$ are continuous 
on $U_{i,\bbC}$ 
with respect to the metric topology.
%
Thus any point $x\in X$ has an open 
neighborhood $U_x$ (in the metric topology) such 
that $U_x\subset U_{i_x,\bbC}$ for some $i_x$, and
\[ \Bigl| \frac{df}{f}(p, v) \Bigr|\le
K_x \max_{j=1,\dots,q_{i_x}}
\Bigl|\frac{d\beta_{i_x,j}}{\beta_{i_x,j}}(p,v)\Bigr| \]
whenever $v \in T_p(X)$ and $p$ is a smooth point of $U_x$ which does
not lie on the
divisors of $f, \beta_{i_x,1}, \dots, \beta_{i_x,q_{i_x}}$. 
The open sets $U_x$ form a cover of $X$; since $X$ compact in 
the metric topology, we can choose a finite subcover
$U_{x_1}, \dots, U_{x_r}$. Now if $K > K_{x_1}, \dots, K_{x_r}$ then
\[ \Bigl| \frac{df}{f}(p, v) \Bigr|\le
K \max_{i,j}\,\Bigl|\frac{d\beta_{ij}}{\beta_{ij}}(p,v)\Bigr| \; . \]
This shows that Proposition~\ref{prop2a} holds.
\qed
\end{reduction}

\begin{reduction} \label{red4} 
%
{\em Suppose $X$ and $X'$ are birationally isomorphic $G$-varieties
over $k\subset\bbC$. 
If Proposition~\ref{prop2a} holds for $X$ and $f\in k(X)$
then it holds for $X'$ and the same $f\in k(X')=k(X)$.}

\smallskip
Indeed, $X$ and $X'$ have isomorphic Zariski-open subsets
$U$ and $U'$. After passing to smaller subsets if necessary, we may 
assume that $U$ and $U'$ are smooth and do not intersect the divisors
of $f, \beta_1, \dots, \beta_n$ on $X$ and $X'$ respectively.
Thus if inequality~\eqref{eqn101} 
holds for every $(p, v)$ such that $p \in U_\bbC$ and $v \in T_p(X_\bbC)$ then
it holds
for every $(p, v)$ such that
$p \in U'_\bbC$ and $v \in T_{p'}(X'_\bbC)$.
The subset $U'_\bbC$ is dense in $X'_\bbC$ with respect to metric topology;
hence, by continuity the same inequality (with the same $\beta_j$ and 
the same $K$) holds for every $(p, v)\in T(X'_\bbC)$ such that $p$
is a smooth point of $X'_\bbC$ and does not
lie in the union of divisors of $f, \beta_1, \dots, \beta_m$. This
means that Proposition~\ref{prop2a} holds for $X'$ as claimed.
\qed
\end{reduction}

We have thus shown that Theorem~\ref{thm2} is a consequence of
the following, purely algebraic statement (see Reductions~\ref{red2}, 
\ref{red3} and~\ref{red4}).

\begin{prop} \label{prop2b} Let $G$ be a finite group, $X$ a projective
$G$-variety, and $f \in k(X)$. Then there exists a 
birational morphism $\pi \colon X' \to X$ of $G$-varieties 
such that the pair $(X', \pi^*(f))$ has property
(*) (see Definition~\ref{def.prop-star}).
\end{prop}

A proof of Proposition~\ref{prop2b} (and thus of Theorem~\ref{thm2})
will be given in the next section.
The idea is to construct $\pi\colon X'\to X$ by resolving
the $G$-action on $X$ to ``standard form'' with respect to a divisor
containing the divisor of $f$; see below.
The simplest (affine) example of such $X'$
is $X'=\bbA^1=\Spec k[t]$,
where $k$ contains a primitive $m$th root of unity, 
$G=\bbZ/m\bbZ$ acts on $\bbA^1$ linearly by a faithful
character, and $f=t$.
In this case we can take $\beta=t^m\in k(X)^G$; the equality
$df/f=\frac{1}{m}d\beta/\beta$
shows that $(X',f)$ has property (*).

\section{Conclusion of the proof}
\label{sect3}

\subsection*{$G$-varieties in standard form}

\begin{defn} \label{def3.1} (\cite[Definition~3.1]{ry})
We say that a generically free $G$-variety $X$ is {\em in standard form
with respect to a divisor} $Y$ if
\begin{itemize}
\item[(a)] $X$ is smooth and $Y$ is a normal crossing divisor on $X$,
\item[(b)] the $G$-action on $X - Y$ is free, and 
\item[(c)] for every $g \in G$ and for every irreducible component $Y_0$ of $Y$
either $g(Y_0) = Y_0$ or $g(Y_0) \cap Y_0 = \emptyset$.
\end{itemize}
\end{defn}

\begin{thm} \label{thm.b1-2} (\cite[Corollary~3.6]{ry})
Let $X$ be a $G$-variety
and $Y\subset X$ be a Zariski closed $G$-invariant
subvariety such that the action of $G$ on $X-Y$ is free. Then
there is a sequence of blowups
\begin{equation} \label{tower1}
\pi\colon X_n \stackrel{\pi_n}{\lra} X_{n-1}
\dots \stackrel{\pi_2}{\lra} X_1 \stackrel{\pi_1}{\lra} X_0 = X 
\end{equation}
with smooth $G$-invariant centers $C_i\subset X_{i}$ such that 
$X_n$ is in standard form with respect to a divisor $\tilde Y$
containing $\pi^{-1}(Y)$.
\end{thm}

\begin{thm} \label{thm:stand.form.properties}
Let $X$ be a $G$-variety in standard form with respect to a divisor
$Y$, let $x$ be a point of $X$, let $Y_1,\dots,Y_m$ be the
irreducible components of $Y$ passing though $x$, and let
$W=Y_1\cap\dots\cap Y_m$. Then
\begin{itemize}
\item[(a)] (\cite[Theorem~4.1]{ry})
$\Stab(x)$ is commutative.
\item[(b)] (\cite[Remark~4.4]{ry})
The action of $\Stab(x)$ on the normal space to $W$ at $x$ is faithful and
decomposes into the sum of one-dimensional representations as follows:
\begin{equation} \label{eqn.rem.31}
T_x(X)/T_x(W)=\bigoplus_{i=1}^m
\frac{T_x(Y_1)\cap\dots\cap\widehat{T_x(Y_i)}\cap\dots\cap T_x(Y_m)}{T_x(W)}\ .
\end{equation}
\item[(c)] (\cite[Remark~4.5]{ry})
Let $e$ be the exponent of $\Stab(x)$.  Then the residue field $k(x)$
of $x$ contains a primitive $e$th root of unity.
\end{itemize}
\end{thm}

\begin{remark} \label{rem:at.point}
Under the assumptions of Theorem~\ref{thm:stand.form.properties}, set
$H=\nvstab(x)$ and $H'=\Stab(x)$. Recall that $H$ acts on $T_x(X)$ 
semi-linearly; see Definition~\ref{def.skew-lin}.
Property (c) of Definition~\ref{def3.1} implies that $Y_i$ is
preserved by the action of $H$, and hence, the subspace
$T_x(Y_i)$ is $H$-invariant for each $i$.  It follows that 
all spaces appearing in \eqref{eqn.rem.31}, are $(k(x)*H)$-modules.

The conormal space $(T_x(X)/T_x(Y_i))^*$ is dual to the $i$th summand
in \eqref{eqn.rem.31}; it is a $(k(x)*H)$-module of dimension $1$ over
$k(x)$, and $H'$ acts on it by a character. Denote this character by
$\xi_i$.  Theorem~\ref{thm:stand.form.properties}(b)
implies that the characters $\xi_1,\dots,\xi_m$ generate the dual 
group $(H')^*$. Combining this observation with
Theorem~\ref{thm:stand.form.properties}(c), we conclude that
Proposition~\ref{prop:one-dim} applies in this setting.
\end{remark}

\subsection*{A local coordinate system}
Suppose that $X$ is an algebraic variety and $x$ is a point of $X$.
Recall that $u_1, \dots, u_n \in \fm_x$ are said to form a local
coordinate system on $X$ at $x$ if their classes 
modulo $\fm_x^2$ form a basis of $\fm_x/\fm_x^2$ 
as a $k(x)$-vector space.

\begin{prop} \label{prop:coord.sys}
Let $X$ be a quasiprojective $G$-variety in standard form with respect
to a divisor $Y$. 
Suppose $Y_1, \dots, Y_m$ are the irreducible components of $Y$ passing though 
a point $x$ of $X$.  Then there exists a local coordinate system
$u_1,\dots,u_m,v_1,\dots,v_l$ at $x$ with the following properties:
\begin{itemize}
\item[(a)] Let $\overline u_i=u_i\bmod\fm_x^2$ and
$\overline v_j=v_j\bmod\fm_x^2$.  Then each of
$\overline u_1,\dots,\overline u_m,\overline v_1,\dots,\overline v_l$
generates a one-dimensional $H$-invariant $k(x)$-subspace of
$\fm_x/\fm_x^2$.
\item[(b)] For every $i = 1, \dots, m$ and every $g \in G$,  
$u_i$ is a local equation of $g(Y_i)$ at $gx$.
\item[(c)] For every $j = 1, \dots, l$ there are integers $e_{j1}, 
\dots, e_{jm} \geq 0$ such that 
$h_j= u_1^{e_{j1}}\dots u_m^{e_{jm}} v_j$ is a
$G$-invariant rational function on $X$.
\end{itemize}
\end{prop}

\begin{pf}
We begin by constructing $u_1, \dots, u_m$.  Consider 
the divisor $D_i =\sum g(Y_i)$, where each summand of the form $g(Y_i)$ 
(for some $g \in G$) appears in this sum exactly once.
Since $X$ is quasiprojective, the divisor $D_i$ can be ``moved off''
the finite set
$Gx$; see~\cite[Theorem~III.1.1]{shaf}.
In other words, for every $i = 1, \dots, m$ there is a rational 
function $u_i$ on $X$ such that the support of the divisor $D_i-(u_i)$ 
does not intersect $Gx$. It is now easy to see that $u_1, \dots, u_m$
satisfy (a) and (b).

Next we turn to the construction of $v_1, \dots, v_l$.
Each $\overline{u}_i$ generates a one-dimensional $k(x)$-subspace
$\lf<\overline u_i\r>=(T_x(X)/T_x(Y_i))^*\subset\fm_x/\fm_x^2$.
We have seen in Remark~\ref{rem:at.point} that $\lf< \overline{u}_i \r>$ 
is $H$-invariant; $H'$ acts on it by the character $\xi_i$.
In view of Lemma~\ref{lem.skew-gr1} and
Proposition~\ref{prop:one-dim}, we can write
\[ \fm_x/\fm_x^2 =  
\lf<\overline u_1\r> \oplus \dots \oplus \lf<\overline u_m\r>
\oplus V_1 \oplus \dots \oplus V_l  \; , \]
where each $V_i$ is a simple $(k(x)*H)$-module and 
$\dim_{k(x)}(V_i) = 1$. 
Choose $\overline v_1,\dots,\overline v_l\in\fm_x/\fm_x^2$ 
so that $\overline v_i$ generates $V_i$ as a $k(x)$-vector space.
Denote the character of $H'$ associated to
$\overline{v}_j$ by $\eta_j$.
By Remark~\ref{rem:at.point}, 
$\xi_1,\dots,\xi_m$ generate the dual group $(H')^*$;
consequently, each $\eta_j$ can be written in the form
\[
\eta_j=\xi_1^{-e_{j1}}\dots\xi_m^{-e_{jm}}
\]
for some integers $e_{j1}, \dots, e_{jm} \geq 0$. 
Note that $\overline h_j=
\overline{u}_1^{e_{j1}}\dots \overline{u}_m^{e_{jm}}\overline{v}_j$
is an $H'$-invariant element of $\fm_x^{d_j}/\fm_x^{d_j+1}$,
where $d_j=e_{j1}+\dots+e_{jm}+1$.
Clearly, $\overline h_j$ generates an $H$-invariant one-dimensional
$k(x)$-subspace $\lf<\overline h_j\r>\subset\fm_x^{d_j}/\fm_x^{d_j+1}$
on which $H'$ acts trivially. By the uniqueness statement in
Proposition~\ref{prop:one-dim}(a), $\lf<\overline h_j\r>\isomo k(x)$
as  $(k(x)*H)$-modules, and by Remark~\ref{rem:triv}, 
after replacing $\overline{v}_j$ 
by $\lambda\overline{v}_j$ for some $\lambda\in k(x)$, we may assume
\begin{equation} \label{e.h-bar}
\overline h_j=
\overline{u}_1^{e_{j1}}\dots \overline{u}_m^{e_{jm}}\overline{v}_j
\text{ \ is an $H$-invariant element of }
\fm_x^{d_j}/\fm_x^{d_j+1}\ .
\end{equation}

We claim that we can choose $v_1, \dots, v_l \in \fm_x$ so that 
$\overline{v}_j = v_j\bmod \fm_x^2$ and each 
$h_j = u_1^{e_{j1}}\dots u_m^{e_{jm}}v_j $ is $G$-invariant.
If we can do this, then $u_1, \dots, u_m, v_1, \dots, v_l$ will clearly
satisfy the requirements of the proposition.

To prove the claim, let $R = \bigcap_{y \in Gx} {\cal O}_{y, X}$ 
be the ring of rational functions on $X$ that are well-defined 
on $Gx$ and let $\cI_{Gx}=\bigcap_{y \in Gx} \fm_y$ be
a $G$-invariant ideal in $R$ consisting of
all elements that vanish on $Gx$.  By our choice of $u_1, \dots, u_m$, 
\begin{equation} \label{e.coordinate}
g^{\ast}(u_i)/u_i \text{ is defined and invertible at every
point of $Gx$}
\end{equation}
for every $g \in G$, and every $i=1, \dots, m$.
Consequently, $u_1^{e_{j1}} \dots u_m^{e_{jm}}\cI_{Gx}$ is 
a $G$-invariant ideal of $R$.

For each $y\in Gx$, the functions $u_1,\dots,u_m$ vanish at $y$, and
hence,
$u_1^{e_{j1}} \dots u_m^{e_{jm}}\cI_{Gx}\subset\fm_{y}^{d_j}$.
Consider the $G$-equivariant projection map 
\[  \textstyle
\psi \colon u_1^{e_{j1}} \dots u_m^{e_{jm}} \cI_{Gx}
\lra \bigoplus_{y \in Gx} \fm_{y}^{d_j}/\fm_{y}^{d_j+1} \; .
\]

If $v\in\cI_{Gx}$ then
$\psi(u_1^{e_{j1}} \dots u_m^{e_{jm}}v)$ depends only 
on the image of $v$ in $\bigoplus_{y \in Gx}\fm_y/\fm_y^2$.
Moreover, since $X$ is 
quasiprojective, the finite set $Gx$ lies in an affine open subset 
of $X$ and hence, by the Chinese Remainder Theorem, the projection map
$\cI_{Gx} \lra \bigoplus_{y \in Gx} \fm_y/\fm_y^2$ is surjective.
Thus
\[ \textstyle
\Im(\psi)=
\bigoplus_{y \in Gx}
\left(u_1^{e_{j1}}\dots u_m^{e_{jm}}\fm_y\right)
\bmod\fm_y^{d_j+1}
\subset\bigoplus_{y \in Gx} \fm_{y}^{d_j}/\fm_{y}^{d_j+1}\ .
\]
We shall denote elements of $\Im(\psi)$ by 
$a =\bigl(a_y\mid y\in Gx\bigr)$, where
$a_y \in  \left(u_1^{e_{j1}}\dots u_m^{e_{jm}}\fm_y\right)
\bmod\fm_y^{d_j+1}$. Recall that by \eqref{e.h-bar}, $\overline h_j$ 
is a nonzero $H$-invariant element of 
$\left(u_1^{e_{j1}}\dots u_m^{e_{jm}}\fm_x\right)
\bmod\fm_x^{d_j+1}$.
Let $a_j$ be the element of $\Im(\psi)$ such that
$(a_j)_y = (g^{-1})^*(\overline h_j)$, where $y = gx$;
in view of \eqref{e.coordinate}, 
\[
(a_j)_y\in
(g^{-1})^*[
\left(u_1^{e_{j1}}\dots u_m^{e_{jm}}\fm_x\right)
\bmod\fm_x^{d_j+1}] = 
\left(u_1^{e_{j1}}\dots u_m^{e_{jm}}\fm_y\right)
\bmod\fm_y^{d_j+1}\ .
\]
Note that since $\overline h_j$ is $H$-invariant, $(a_j)_y$ is independent
of the choice of $g$.
By our construction $a_j$ is $G$-invariant and $(a_j)_x = 
\overline h_j \in \fm_x^{d_j}/\fm_{x}^{d_j+1}$.

The homomorphism $\psi$ has a $G$-equivariant $k$-linear
splitting and consequently, there exists a $G$-invariant element 
$h_j = u_1^{e_{j1}}\dots u_m^{e_{jm}}v_j \in
 u_1^{e_{j1}}\dots u_m^{e_{jm}}\cI_{Gx}$ such that
$\psi(h_j) = a_j$.
In particular, 
$\overline h_j=
\overline{u}_1^{e_{j1}}\dots\overline{u}_m^{e_{jm}}
\overline{v}_j= h_j \bmod{\fm_x^{d_j+1}}$ and hence 
$\overline{v}_j=v_j\bmod \fm_x^2$. This proves the claim and thus
shows that $u_1, \dots, u_m, v_1, \dots, v_l$ have the required properties.
\end{pf}

\subsection*{Property (*)}

We are now ready to 
revisit property (*) of Definition~\ref{def.prop-star}.

\begin{lem} \label{lem:diff.poles}
Suppose $X$ be a quasiprojective $G$-variety in standard form with respect
to a divisor $Y$, $x\in X$, $u_1, \dots, u_m$, 
$v_1, \dots, v_l$, and $h_1, \dots, h_l$, are as in
Proposition~\ref{prop:coord.sys}, and
$w_i=\prod_{g\in G}g^*(u_i)$. Let $f$ be a rational function
on $X$ whose divisor is supported on $Y$. Then 
\[  
df/f \in\cM \; , \] 
where $\cM$ is the $\cO_{x,X}$-module generated by
$dw_i/w_i$ and $dh_j/h_j$ with $i=1,\dots,m$ and $j=1,\dots,l$.
\end{lem}

Note that $w_i$ and $h_j$ are $G$-invariant rational 
functions on $X$ for every $i=1, \dots, m$ and $j = 1, \dots, l$.

\begin{pf}
Let $(\Omega_X^1)_x$ be the $\cO_{x,X}$-module of germs at $x$ of
regular differential forms on $X$.
Since $X$ is smooth,
$(\Omega_X^1)_x$ is a free $\cO_{x,X}$-module generated by
$du_1,\dots,du_m,dv_1,\dots,dv_l$.  Let $\cM'$ be the $\cO_{x,X}$-module 
generated by $du_i/u_i$, and $dv_j/v_j$, 
where $i=1,\dots,m$ and $j=1,\dots,l$.  Clearly,
$(\Omega_X^1)_x\subset\fm_x\cM'$.

We claim that $\cM=\cM'$.
It is clear that $\cM\subset\cM'$;
to prove the opposite inclusion,
we shall show that $dw_1/w_1,\dots,dw_m/w_m,dh_1/h_1,\dots,dh_l/h_l$
generate $\cM'$ as an  $\cO_{x,X}$-module.

Note that by our choice of $u_1, \dots, u_m$, we  can write
$w_i=a_iu_i^{|G|}$ for some $a_i \in {\cal O}_{X, x} - \fm_x$; 
see~\eqref{e.coordinate}.  Thus
\[
\frac{dw_i}{w_i} = \frac{da_i}{a_i} + |G| \frac{du_i}{u_i}
\equiv|G|\,\frac{du_i}{u_i}\pmod{(\Omega_X^1)_x}\ .
\]
In particular, since $(\Omega_X^1)_x\subset\fm_x\cM'$, we conclude that
\begin{equation} \label{e.dw_i}
\frac{dw_i}{w_i} 
\equiv|G|\,\frac{du_i}{u_i}\pmod{\fm_x\cM'}\ .
\end{equation}
On the other hand, since $h_j=u_1^{e_{j1}} \dots u_m^{e_{jm}} v_j$,
we have
\begin{equation} \label{e.dh_j}
\frac{dh_j}{h_j}= \frac{dv_j}{v_j}+\sum_ie_{ji} \frac{du_i}{u_i}\ .
\end{equation}
Examining~\eqref{e.dw_i} and~\eqref{e.dh_j}, we see
that $dw_i/w_i$ ($i=1,\dots,m$), and $dh_j/h_j$
($j=1,\dots,l$) generate $\cM'/\fm_x\cM'$ as a $k(x)$-vector space.
Consequently, by Nakayama's lemma these elements generate $\cM'$
an $\cO_{x,X}$-module.
Thus $\cM'=\cM$, as claimed.

Since the divisor of $f$ is supported on $Y$, locally near $x$ 
it is a union of smooth hypersurfaces of the form $\{u_i=0\}$. 
This means that $f = au_1^{e_1} \dots u_m^{e_m}$ for some 
$a \in {\cal O}_{x, X} - \fm_x$ and $e_1, \dots, e_m \geq 0$; hence,
\[ \frac{df}{f} = \frac{da}{a} + e_1\frac{du_1}{u_1} + \dots 
+ e_m \frac{du_m}{u_m} \in\cM'=\cM \; . \]
\end{pf}

\begin{prop} \label{prop.star} 
Let $X$ be a projective variety in standard form with respect 
to a divisor $Y$ and let $f$ be a rational function on $X$ whose divisor 
is supported on $Y$. Then the pair $(X, f)$ has property (*); see 
Definition~\ref{def.prop-star}.
\end{prop}
 
\begin{pf}
By Zariski compactness, it is enough to show that for any $x\in X$, 
there exist $\beta_{1}, \dots, \beta_{q} \in k(X)^G$ and
$\gamma_{1},\dots,\gamma_{q}\in \cO_{x,X}$ such that
\[ \frac{df}{f} = \gamma_{1}\frac{d\beta_{1}}{\beta_{1}} + \dots
\gamma_{q}\frac{d\beta_q}{\beta_q}  \; . \]
The last assertion is immediate from Lemma~\ref{lem:diff.poles}.
\end{pf}

\subsection*{Proof of Theorem~\ref{thm2}}

As we showed in Section~\ref{sect2}, 
it is enough to prove Proposition~\ref{prop2b}.

Assume $X$ be a projective $G$-variety and $f \in k(X)$. Let
$X_0$ be the subvariety of all points in $X$ with nontrivial
stabilizers.  Let $Y$ be the union of  $X_0$ and (the supports of)
the divisors of $g^{\ast}(f)$ for every $g \in G$; it is a
$G$-invariant Zariski closed subvariety of $X$.
By Theorem~\ref{thm.b1-2} there exists a birational 
morphism $\pi \colon X' \lra X$ and a divisor $Y'\subset Y$,  
such that $X'$ is in standard form with respect to $Y'$ and
$\pi^{-1}(Y) \subset Y'$. Note that
the divisor of $\pi^{*}(f)$ is contained in $\pi^{-1}(Y)$ 
and hence in $Y'$.
Proposition~\ref{prop2b} (and thus Theorem~\ref{thm2}) now follows
from Proposition~\ref{prop.star}, which asserts that the pair 
$(X', \pi^{\ast}(f))$ has property (*).
\qed

\end{document}